\documentclass[11pt]{article}
\usepackage[usenames]{color}
\usepackage[leqno]{amsmath}
\usepackage{amssymb,amsfonts,amstext,amsthm,amsmath}
\newcommand{\R}{\mathbb{R}}\newcommand{\ov}{\overline}
\newcommand{\nb}{\nabla}\newcommand{\pa}{\partial}

\newcommand{\inp}[2]{\left\langle #1,\!#2\right\rangle}
\renewcommand{\div}{{\rm div}}
\renewcommand{\d}{{\rm d}}
\renewcommand{\cal}{\mathcal}
\newcommand{\tr}{{\rm c}}

\renewcommand{\u}{\nb u}

\newcommand{\GV}{\nb V}
\newcommand{\GVS}{(\nb V)^{\ast}}
\newcommand{\Ric}{{\rm Ric}}
\newcommand{\gz}{\gamma_0}\newcommand{\gm}{\gamma}
\newcommand{\sg}{\sigma}\newcommand{\bsg}{\overline{\sigma}}
\newcommand{\M}{\mathcal{M}}
\newcommand{\al}{\alpha}\newcommand{\be}{\beta}
\newcommand{\bvab}{BV(\M,\{\al,\be\})}
\newcommand{\bvabs}{BV_s(\M,\{\al,\be\})}
\newcommand{\Ez}{\mathcal{E}_0}\newcommand{\Ee}{\mathcal{E}}
\newcommand{\brho}{\bar{\rho}}
\newcommand{\vv}{\mathcal{V}}\newcommand{\vz}{\mathcal{V}_0}
\newtheorem{thm}{Theorem}
\newtheorem{lem}[thm]{Lemma}
\newtheorem{cor}[thm]{Corollary}
\theoremstyle{definition}
\newtheorem*{dfn}{Definition}
\theoremstyle{remark}
\newtheorem{rem}{Remark}

\newcommand{\N}{\cal{N}}
\newcommand{\tV}{\tilde{V}}

\author{Arnaldo S. Nascimento(1) and  Alexandre C. Gon\c{c}alves(2)
\thanks{Partially supported by FAPESP grant \#2006/02023-4}}

\title{ Instability results for an elliptic equation
on compact Riemannian manifolds with non-negative Ricci curvature}

\begin{document}

\maketitle

\begin{abstract}
 We prove nonexistence of nonconstant local minimizers for a  class of
functionals, which typically  appears  in the scalar two-phase field
model, over a smooth $N-$dimensional Riemannian manifold without
boundary with non-negative Ricci curvature. Conversely for a class
of surfaces possessing a simple closed geodesic along which the
Gauss curvature is negative we prove existence of nonconstant local
minimizers for the same class of functionals.
\end{abstract}



\section{Introduction}

Let $\M$ be a smooth $N-$dimensional compact Riemannian manifold without
boundary and consider the  functional $\Ee:H^1(\mathcal{M})\to\R$
given by
   \begin{equation}\label{funcional-E}
    \mathcal{E}(u)=\int_{\mathcal{M}}\{ \frac{|\nabla
    u|^2}{2}-F(u)\}\,d\mu
    \end{equation}
where $F$ is a $C^2$ real function and $H^1(\mathcal{M})$ the usual
Sobolev space.

In this work we are interested in the question of how locally
minimizing functions of $\Ee$ are related to the geometry of
$\M$.

We will say that $u_0 \in C^\infty(\mathcal{M})$ is a local
minimizer of $\mathcal{E}$ if $\exists \, \delta>0$ such that $$
\mathcal{E}(u_0) \leq \mathcal{E}(u)\;\;\;\mbox{whenever}\;\;
\|u-u_0\|_{H^1(\mathcal{M})}\leq \delta. $$

In case the first inequality  is strict, i.e., $\Ee(u_0) < \Ee(u)$,
$u_0$ is said to be a local isolated minimizer.
Our main results are stated in the following
theorems.
\begin{thm}\label{tunstable}
Suppose that the Ricci curvature of $\M$ is non-negative.
Then any local minimizer of $\Ee$ is a constant function.
\end{thm}

An interesting condition that shows up in the computations of
Theorem \ref{tunstable} provides some insight
 on the structure of $\M$. For any $u\in
H^1(\M)$ we denote by $\Ee''(u)$ the second variation of $\Ee$ at
$u$.
\begin{thm}\label{tsemistable}
Keep the hypothesis of Theorem \ref{tunstable}. Let $u$ be a
non-constant critical point of $\Ee$ and set $v=|\u|$. If
\begin{equation}
   (\Ee''(u)v,v)=0
\end{equation}
then there exist a complete riemannian submanifold $\N\subset\M$, a
real geodesic line subundle $\cal{I}\subset T\M$ and an isometric
regular covering map $\varphi:\R\times\N\to\M$. Denoting by $K$ the
group of covering transformations of $\varphi$, then $K$ is made of
isometries, and $\M$ is isometric to the quotient $(\R\times\N)/K$.
If $\cal{I}$ is orientable then $K$ is generated by a nontrivial
(affine) translation of $\R$ with some isometry of
$\N$. Otherwise $K$ is generated by two involutions of
$\R\times\N$.
\end{thm}

Regarding Theorem \ref{tunstable} we show how to construct non-constant local
minimizers on some non-convex surfaces. To that purpose we introduce
a small positive parameter $\varepsilon$ in the functional  thus
writing
\begin{equation}\label{variation-eps}
\Ee_\varepsilon(u)=\int_{\M}\{\varepsilon
\frac{|\nabla u|^2}{2}-\varepsilon^{-1} F(u)\}\,d\mu
\end{equation}
and  take for $F$ a suitable nonnegative double-well potential which
vanishes only at $\alpha$ and $\beta\,\, (\alpha<\beta)$. As usual
$\chi_A$ will stand for the characteristic function of a set $A$.

\begin{thm}\label{tstable}
Let  $\M$ be a surface diffeomorphic to $S^2$. Assume that there
exists a simple closed geodesic $\gz\subset\M$ so that the Gauss
curvature $K$ of $\M$ is negative along $\gz$. Then for
$\varepsilon$ small enough there is a non-constant family
$\{u_\varepsilon\}_{\varepsilon >0}$ of local minimizers of
$\Ee_\varepsilon$. Moreover it holds that $u_\varepsilon
\stackrel{\varepsilon \to 0}{\longrightarrow}u_0$ in $L^1(\M)$ where
$u_0=\al\chi_{\M_{\al}}+\be\chi_{\M_{\be}}$ and
$\M=\M_{\al}\cup\gz\cup\M_{\be}$ is the partition of $\M$ determined
by $\gz$ .
\end{thm}
The association of local minimizers  of $\Ee$ to the
geometry of the domain goes back to 1978 when the authors in
\cite{CH} and \cite{M} considered the evolution problem
\begin{equation}\label{CH-M-eq}
   \left\{\begin{array}{l}
       \displaystyle
 u_t=\triangle u \!+\!f(u)\;\;\;\; \mbox{in}\ \ \R^+\times \Omega\
\\
         \displaystyle
 \partial_{\nu}u = 0\;\;\;\;\;\;\;\;\;\;\;\;\;\;\; \mbox{on}\ \
\R^+\times \partial\Omega\,
\end{array}\right.
\end{equation}
where $\Omega \subset \R^N$ is a smooth bounded domain,  $f \in
C^2(\Omega)$ and $\partial_{\nu}$ stands for the exterior normal
derivative.

They showed that if $\Omega$ is convex then any non-constant
solution to \eqref{CH-M-eq} is unstable in the Lyapunov sense. In
this case it amounts to saying that any local minimizer of the
corresponding energy functional is a constant function.

Still for bounded convex domains with homogeneous zero Neumann
boundary condition, the same kind of result was obtained for
systems of reaction-diffusion equations \cite{JM} and \cite{Lopes},
Ginzburg-Landau equation \cite{JS}, reaction-diffusion systems with
skew-gradient structure \cite{Y}, geometric parabolic equation
\cite{JZ} and in the context of permanent currents for the full
bi-dimensional  Ginzburg-Landau functional in \cite{JS}, among
others. In all of these works the  proofs make use in a strong way
of the homogeneous Neumann boundary condition on a convex domain.

When $\M$ is a general Riemannian boundaryless manifold the
Euler-Lagrange equation for $\Ee$  yields the stationary solutions
of the reaction-diffusion equation
\begin{equation}\label{ut-parab}
    \displaystyle
 u_t=\Delta u \!+\!f(u)\;\;\;\; \mbox{in}\ \ \R\times\M\ .
\end{equation}
The only result of this type regarding \eqref{ut-parab} over surfaces was
considered in \cite{RW} where it was shown that if $\M\subset\R^3$
is a convex surface of revolution then the only stable solutions are
the constant ones. Actually the prove consists of showing that
\eqref{funcional-E}, with $F'=f$, has no nonconstant local
minimizer.

In this particular case writing the planar curve that generates the surface in
appropriate coordinates reduces the domain to an interval thus
making the underlying analysis much easier than the general case
considered here.

In case $\mathcal{M}$ is a bounded domain in $\R^N$ typically
$\Ee_\varepsilon$  models the phase separation phenomenon in the
context of van der Waals-Cahn-Hilliard theory whereby $u$ represents
the density of a two-phase fluid and is also associated to
the motion of phase boundaries (interfaces) by mean curvature
(see \cite{I}, for instance).

Equation \eqref{ut-parab} has been studied in the context of pattern formation,
i.e., existence of nonconstant stable (in the sense of Lyapunov)
stationary solution. It may model bio-chemical processes over cell
surfaces or propagation of calcium waves over the surface of a
fertilized egg, for instance.

In particular Theorem \ref{tunstable} implies that \eqref{ut-parab}
has no pattern as long as $\M$ has non-negative Ricci curvature.
On the other hand Theorem \ref{tstable} gives an example of $\M$ for which
\eqref{ut-parab}, after a suitable scaling, develops patterns.

Setting $f=F'$ then clearly critical points of $\Ee$ satisfy the
semi-linear elliptic equation
\begin{equation}\label{autoeq}
   \Delta u + f(u)=0\;\mbox{on}\;\M.
\end{equation}
A smooth solution $u$ of the above equation is said to be weakly stable
if the quadratic form
\begin{equation}
 E(\varphi)=\int_{\mathcal{M}}\{ \frac{|\nabla
      \varphi|^2}{2}-f'(u)\varphi^2\}\,d\mu \geq 0,
\end{equation}
in $H^1(\M).$  Otherwise $u$ is called weakly unstable.
Then it follows immediately from the proof of Theorem \ref{tunstable}
that any nonconstant solution to the above equation is weakly unstable
as long as $\M$ has non-negative Ricci curvature.

The paper is organized as follows. In Section 2 in addition to
recalling  some notation of Riemmanian Geometry we prove some
preliminary results, Section 3 is  devoted to the proofs of Theorem
\ref{tunstable} and Theorem \ref{tsemistable} and Section 4 to
the proof of Theorem \ref{tstable}.

\section{Geometric Background}

Throughout this section $\M$ will denote an $N$-dimension ($N \geq
2$) riemannian manifold without boundary, and $T\M$, $T^\ast\M$
its tangent and cotangent bundles, respectively. We shall deal
with the tensor bundles $T^r_s(\M)=(T\M)^{\otimes r}\otimes (T^\ast
\M)^{\otimes s}$, for  non-negative integers $r$ and $s$. For an
integer $k\geq 0$ let $\cal{A}^kT^\ast\M$ be the alternate k-bundle
of $T^\ast\M$. Notice that $T^0_0(\M)=\cal{A}^0T^\ast\M$ is the
trivial bundle $\M\times\R$.

Given  any real vector bundle $\cal{F}$ over $\M$ we denote by
$G(\cal{F})$ the set of its smooth sections and by
$G^k(\cal{F})=G(\cal{A}^kT^\ast\M\otimes\cal{F})$ the smooth
sections of \emph{k-forms of $\M$ with coeficients in $\cal{F}$}.

The contraction is a natural coupling $c:T^1_1(\M)\to
T^0_0(\M)$ given by $c(v\otimes\omega)=\omega(v)$, where
$v\otimes\omega$ is a decomposable tensor of $T\M\otimes
T^\ast \M$. The contraction extends to
$c:T^r_s(\M)\to T^{r-1}_{s-1}(\M)$ for any
$r,s\geq 1$, by putting $c(v_1\otimes\dots\otimes
v_r\otimes\omega_s\otimes\dots\otimes \omega_1)=\omega_1(v_1)\
v_2\otimes\dots\otimes v_r\otimes\omega_s\otimes
\dots\otimes\omega_2$. Indeed, when $r=s=1$ the contraction is
just the trace operator on linear homomorphisms $T\M\to
T\M$.

Let $\nb:G(T^1_0(\M))\longrightarrow G^1(T^1_0(\M))$ be
the Levi-Civita connection on $\M$.
It is well known that $\nb$ can be extended in a
unique way to an operator $\ov{\nb}:G(T^r_s(\M))\longrightarrow
G^1(T^r_s(\M))$ such that Leibnitz rule is preserved
and commutes with the contraction (\cite{KN}).
We abuse notation and write
$\overline{\nb}=\nb$ whenever $r,s$ are not both zero.
When $f\in G(T^0_0(\M))$ is just a smooth function we preserve the
usual notation $\nb f=(df)^\ast\in G(T^1_0(\M))$. It then
follows
\begin{align}
 \nb(T\otimes W)&=\nb T\otimes W+
 T\otimes\nb W  \label{tensor-1}\\
     \forall\ T&\in G(T^r_s(\M))\ \text{and}\ \forall\ W\in G(T^p_q(\M))\ ,
\ \text{and}\notag\\ \nb c(T)&=c(\nb T)\ ,\label{tensor-2}\\
 &\text{for a contraction}\ c:T^r_s(\M)\rightarrow T^{r-1}_{s-1}(\M)\ .
\notag
\end{align}

Notice that we identify
\begin{equation}\begin{aligned}
(T\M)^{\otimes r}\otimes(T^\ast\M)^{\otimes s}&\otimes
(T\M)^{\otimes p}\otimes(T^\ast\M)^{\otimes q}\cong\\
&(T\M)^{\otimes r}\otimes(T\M)^{\otimes p}\otimes
(T^\ast\M)^{\otimes s}\otimes(T^\ast\M)^{\otimes q}\ ,
\end{aligned}\end{equation}
and similarly, by sticking the 1-form component of a section
of $\cal{A}^1T^\ast\M\otimes (T^r_s(\M))$ on the left of the covariant
part we have
$\cal{A}^1T^\ast\M\otimes T^r_s(\M)\cong T^r_{s+1}(\M)$.
These identifications are necessary for \eqref{tensor-1}
and \eqref{tensor-2} to make sense. They also allow us to define
the composition $\nb(\nb T)$ for any $T\in G(T^r_s(\M))$.

Some combinations of $\otimes$ and $c(\cdot)$ deserve special
notation. For tensors $T\in G(T^1_s(\M))$ and $W\in G(T^1_q(\M))$
we write $TW=c(W\otimes T)$. When $s=1$ and $q=1$, $TW$ is
the composition of the endomorphisms $T$ with $W$, and if $q=0$
$TW$ is the image of the vector $W$ under $T$. In particular,
if $s\geq2$ and $W_1,W_2$ are vector fields we set
$T(W_1,W_2)=[TW_2]W_1$.

Let $F\in T^1_3(\M)$ be the Riemann tensor of
$\M$. The tensor $F$ can be seen as a two form with
values in the endomorphism bundle of $T\M$ or $F\in
G^2(T\M\otimes T^\ast\M)$.
For any vector fields $X,Y,Z$ and $W$ locally defined we have
\begin{equation}\begin{aligned}
F(X,Y,Z,W)&\stackrel{{\rm def}}{=}\inp{[FZ](Y,X)}{W}=\\
        &=\inp{\nb_X\nb_YZ-\nb_Y\nb_XZ-\nb_{[X,Y]}Z}{W}\ .
\end{aligned}\end{equation}
The proof of the next lemma is straightforward and will be
omitted.
\begin{lem}
Let $V\in G(T^1_0(\M))$. Then the skew-symmetric
component respect to the cotangent factors of $\nb(\nb V)$ is
$FV$. This is equivalent to
$$[\nb(\GV)](X,Y)-[\nb(\GV)](Y,X)=[FV](Y,X)$$ for any vectors
$X,Y$.
\end{lem}

We define the Ricci tensor of $\M$ as $\Ric(V,W)=-c([FW]V)$,
for any $V,W$ vector fields. Observe that if $\{s_i\,|\,i=1,\dots,n\}$
is any local orthonormal basis of $T\M$ then
$\Ric(V,W)=\sum_{i=1}^n F(s_i,V,W,s_i)$.
\begin{dfn}
A non-negative Ricci manifold $\M$ is
one that satisfies \\ $\Ric(V,V)\geq 0$ for any $V\in T\M$.
\end{dfn}

The following lemma will be useful in our approach.
\begin{lem}\label{lem-anti-sim}
Let $V$ and $W$ be vector fields over $U\subset \M$ open.
Then
\begin{equation}\label{anti-sim}
  \tr([\nb(\GV)]W-\nb_W(\GV))=\Ric(W,V)\ .
\end{equation}
\end{lem}
\begin{proof}[Proof]
We choose an orthonormal basis $\{s_1,s_2,\dots,s_n\}$ locally
defined and compute
\begin{align*}
  \tr&([\nb(\GV)]W-\nb_W(\GV))=\\
   &=\sum_i \inp{[\nb_{s_i}(\GV)]W-[\nb_W(\GV)]s_i}{s_i}=\\
   &=\sum_i \inp{\nb_{s_i}[(\GV)W]-(\GV)\nb_{s_i}W-\nb_W[(\GV)s_i]
+(\GV)\nb_Ws_i}{s_i}=\\
   &=\sum_i \inp{\nb_{s_i}\nb_WV-\nb_W\nb_{s_i}V-\nb_{[s_i,W]}V}{s_i}=\\
   &=\sum_i F(s_i,W,V,s_i)=\Ric(V,W)\ .\qedhere
\end{align*}
\end{proof}

Let $\M$ be a Riemann surface and $\gz\subset\M$ be a simple
closed geodesic. We assume local orientability of $\M$ in a neighborhood
of $\gz$, i.e., there exists a smooth unitary orthogonal vector
field $\eta$ defined on $\gz$. Standard arguments (see \cite{dC})
allow $\eta$ to be extended to a geodesic vector field on a vicinity
$\vv$ of $\gz$. Let $\varphi_t(p)=\varphi(t,p)$ be the
flow of $\eta$. Restricting $\vv$ if necessary, one can choose $\delta>0$
so that the map $\varphi:[-\delta,\delta]\times\gz\to \vv$
is a diffeomorphism. In all computations it is implicitly assumed that 
$\gz$ is arcwise parametrized, so that $\gz'$ is a well defined 
unitary vector field \emph{along} $\gz$.

Let $t$ and $x$ be the coordinate functions of the inverse map
$\varphi^{-1}:\vv\to[-\delta,\delta]\times\gz$, $\varphi^{-1}(p)=(t(p),x(p))$.
For any $\sg:[0,1]\to \vv$ a smooth curve we denote by $\bsg$ its
projection over $\gz$,
\begin{equation}
   \bsg(s)=x\circ \sg(s)\ ,\qquad 0\leq s\leq 1\ .
\end{equation}
Notice that we abuse language and denote by $\sg$ either a curve or
its trace, according to the context. Similarly, $|\sg|$ denotes the
length of the curve, but for a two dimensional region $U\subset\M$,
$|U|$ denotes its area.

The contents of the next lemma are well known to geometers, and can
be found in the literature. Nevertheless we choose to state and
proof the precise statements we need for the sake of completeness.
\begin{lem}\label{length-K}
Suppose that the gaussian curvature $K$ is strictly negative on $\vv$.
We have:\\
({\rm a}) Let $p_0,p_1\in \vv$ and $\sg$ be any smooth simple curve
joining $p_0$ and $p_1$. Then\\
({\rm a1}) $|\sg|\geq|t(p_1)-t(p_0)|$ and equality holds if and only
if $\sg$ reparametrizes the geodesic segment $t\mapsto\varphi_t(p)$ 
between $p_0$ and $p_1$.\\
({\rm a2}) $|\sg|\geq|\bsg|$ and equality holds if and only if
$\sg=\bsg\subset\gz$.\\
({\rm b}) Let $J\subset\gz$ be an interval or $J=\gz$.
Let $0<\delta_0\leq\delta$ and $U$ be any of the sets
$\varphi([0,\delta_0]\times J)$ or $\varphi([-\delta_0,0]\times J)$.
Then $|U|>\delta_0\,|J|$.\\
\end{lem}
\begin{proof}
Let $W_p=(d\varphi_t)_x\cdot\gz'(x)$ for any 
$p=\varphi(t,x)\in \vv$. Then $W$ is a smooth
vector field on $\vv$. Using the symmetry of the Levi-Civita
connection together with $|\eta|\equiv 1$ one gets
\begin{align}
   \frac d{dt}\inp{\eta}{W}&=\eta\inp{\eta}{W}=
 \inp{\nb_{\eta}\eta}{W}+\inp{\eta}{\nb_{\eta}W}=\\
&=\inp{\eta}{\nb_W\eta+[\eta,W]}=\\
&=\frac12W\inp{\eta}{\eta}+\inp{\eta}{[\eta,W]}=0\ ,
\end{align}
and therefore $\inp{\eta}{W}$ is constant along the flow of $\eta$.
Over $\gz$ we know $W=\gz'$, from what we obtain $\inp{\eta}{W}=0$ on $\vv$.
The field $W$ is nowhere singular in $\vv$, and we set the orientation
of $\vv$ as given by the orthogonal basis $\{\eta,W\}$.

Let $x(s)$ and $t(s)$ be the local coordinate functions of a given
$\sg(s)$, so that $\sg(s)=\varphi(t(s),x(s))$, 
for $0\leq s\leq 1$. Let $p_0=\sg(0)$ and $p_1=\sg(1)$.
Notice that $x(s)$ belongs to the trace of $\gz$ and its derivative is 
a multiple of $\gz'$, but we abuse language and set 
$x'(s)=\inp{\frac{dx}{ds}}{\gz'}_{x(s)}$. Then $\sg'(s)=t'(s)\eta+x'(s)W$.
Since $\sg$ has no self-intersections it follows
\begin{equation}\begin{aligned}\label{comp-sg}
 |\sg|=\int_0^1|\sg'(s)|\,ds &= \int_0^1\sqrt{(t')^2+(x')^2|W|^2}\,ds\\
            &\geq\int_0^1|t'|\,ds \geq |t(p_1)-t(p_0)|\ .
\end{aligned}\end{equation}
Equality in (15) occurs if and only if $x'\equiv0$ and $t'$ does not
change sign. This implies $x(s)=x_0\in \gz$ is constant, hence
$\sg(s)=\varphi_{t(s)}(x_0)$ is just a parametrization of an arc of
geodesic, what proves part (a1). A computation similar to
\eqref{comp-sg} yields
\begin{equation}\label{comp-sg-x}
|\sg|=\int_0^1\sqrt{(t')^2+(x')^2|W|^2}\,ds\geq\int_0^1|x'||W|\,ds\ .
\end{equation}
We show that $|W_p|\geq 1$ with equality only when $p\in\gz$. It
suffices showing that the function $t\mapsto |W_{(t,x)}|^2$ is
convex in $[-\delta,\delta]$, with a strict minimum attained in
$t=0$. Indeed,
\begin{equation}
\left.\frac d{dt}\right|_{t=0}|W|^2=2\inp{\nb_{\eta}W}{W}_{t=0}=
   2\inp{\nb_W\eta}{W}_{t=0}=0\ ,
\end{equation}
since $\nb_W\eta=0$ along $\gz$. The second derivative gives
\begin{align}
 \frac{d^2}{dt^2}|W|^2 &=
   2\left(\inp{\nb_{\eta}\nb_W\eta}{W}+\inp{\nb_W\eta}{\nb_{\eta}W}\right)\\
 &=2\left(\inp{\nb_{\eta}\nb_W\eta-\nb_W\nb_{\eta}\eta-\nb_{[\eta,W]}\eta}{W}+
|\nb_W\eta|^2\right)\\
 &=2(-K+|\nb\eta|^2)|W|^2\ ,
\end{align}
 so it is strictly positive for any $t$, under the hypothesis $K<0$.
Back to \eqref{comp-sg-x} we have
\begin{equation}
|\sg|\geq\int_0^1|x'||W|\,ds\geq\int_0^1|x'|\,ds=|\bsg|\ ,
\end{equation}
with equality $|\sg|=|\bsg|$ if and only if $t(s)\equiv0$, or
$\sg=\bsg$ is an arc of the geodesic $\gz$. This proves (a2).

Now assume $U=\varphi([0,\delta_0]\times J)$.
We consider an orthonormal basis of 1-forms $\{\omega_1,\omega_2\}$
dual to $\{\eta,\frac W{|W|}\}$. The area element is $\omega_1\wedge
\omega_2$. Let $J\subset\gz$ be arclength parametrized by the
interval $[s_0,s_1]\subset\R$ so that $\{\eta,J'(s)\}$ preserves
the orientation of $\vv$ over $\gz$. Using the local chart $\varphi$
to write the integral on the plane and applying Fubini Theorem
the area of $U$ is computed as
\begin{equation}\begin{aligned}
  |U|=\int_U\omega_1\wedge\omega_2&=\int_{[0,\, \delta_0]\times[s_0,\, s_1]}
\varphi^{\ast}(\omega_1\wedge\omega_2)=\\
&=\int_{s_0}^{s_1}\int_0^{\delta_0}|W|\,dtds>
    \int_{s_0}^{s_1}\int_0^{\delta_0}\,dtds=\delta_0|J|\ ,
\end{aligned}\end{equation}
thus proving the theorem.
\end{proof}

\section{Nonexistence of nonconstant minimizers }

This section is devoted to the proofs of Theorems \ref{tunstable}
and \ref{tsemistable}, which  in turn will be applications of
the identities established in the next two lemmas.

Recall that the riemannian metric of $\M$ induces metrics
in any tensor product $T^r_s(\M)$, as well as in their spaces
of sections. If $T,W\in T^1_1(\M)$ then their inner-product (fiberwise)
is computed as $\inp{T}{W}=c(c(T\otimes W^\ast))$, being $W^\ast$ the
(metric) transpose of the endomorphism $W:T\M\to T\M$.

If $V$ is a $C^1$ vector field on $\M$ we set $\div(V)=c(\GV)$.
The hessian of a $C^2$ function $u$ on $\M$ is $H_u=\nb(\u)$.
The laplacean of $u$ is then $\Delta u=c(H_u)=\div(\u)$.

The riemannian measure on $\M$ will be denoted by $d\mu$.
By a \emph{component} of a topological space we always mean a
\emph{connected component}.

\begin{lem}\label{ident-V}
Let $V$ be a $C^2$ vector field on $\M$ and $u$ a $C^3$
function on $\M$. Then
\begin{equation}\label{identidade}
  \Delta(Vu)-V(\Delta u)= \div(\GVS\u)+\inp{H_u}{\GV}+\Ric(\u,V)\ .
\end{equation}
\end{lem}
\begin{proof}[Proof]
We first notice that
\begin{equation}\label{formI}
   \nb(Vu)=[\d(Vu)]^\ast=\GVS\u+H_u\,V\ .
\end{equation}
Then,
\begin{align}
  \Delta(Vu)&-V(\Delta u)=\notag\\
     &=\tr(\nb[\GVS\u+H_u\,V])-\tr(\nb_VH_u)\notag\\
     &=\tr(\nb[\GVS\u])+\tr([\nb H_u]V+H_u\GV-\nb_VH_u)\notag\\
     &=\div(\GVS\u)+\tr([\nb H_u]V-\nb_VH_u)+\tr(H_u\GV)\label{ident-pr}\ .
\end{align}
Applying Lemma \ref{lem-anti-sim} to the second summand of term
\eqref{ident-pr} and observing that $\tr(H_u\GV)=\inp{H_u}{\GV}$
we arrive at
\begin{equation*}
   \Delta(Vu)-V(\Delta u)=\div(\GVS\u)+\Ric(\u,V)+\inp{H_u}{\GV}\ ,
\end{equation*}
and the proof is complete.
\end{proof}
\begin{rem} Lemma \ref{ident-V} is central in the next constructions of
this section. Indeed, it somehow appears in \cite{RW}, where
its full geometric significance is shadowed by the high symmetry
of that case. The main idea there, which holds in general,
is a commutation relation between the laplacean operator and a
particular directional derivative, namely, the normalized gradient of $u$.
\end{rem}

Let $u$ be a non-constant critical point  of $\Ee$ with $F'=f$. Then
\begin{equation}
\frac{d}{dt}\Ee(u+t v)|_{t=0}=-\int_\M(\Delta u+f(u))\,v\,d\mu=0,
 \;\forall\;v \in H^1(\M)\ .
\end{equation}
The linearization of the operator $\Delta+f(\cdot)$ at $u$
yields an operator $\cal{L}:H^1(\M)\to H^{-1}(\M)$ defined by
\begin{equation}
   \cal{L}(u)v=\Delta v+i(f'(u)\,v)\ ,
\end{equation}
where $i:H^1(\M)\to H^{-1}(\M)$ is the Sobolev inclusion $H^1\subset H^{-1}$.
Let $(\cdot,\cdot):H^{-1}\times H^1\to\R$ be the canonical pairing
of a vector space and its dual. Then
\begin{equation*}
\frac{d^2}{dt^2}\Ee(u+tv)|_{t=0}=(\Ee''(u)v,v)=-(\cal{L}(u)v,v)\ .
\end{equation*}
For the next lemma we temporarily drop any hypothesis about
Ricci curvature. It will be imediate that for Ricci non-negative
manifolds the quadratic form associated to $\cal{L}$
is not sign definite. Define
\begin{equation*}
  U \stackrel{{\rm def}}{=} \{\u\neq0\} \subset \M\ .
\end{equation*}
Let $V$ be the unitary vector field $V=\frac{\u}{|\u|}$ over $U$.
\begin{lem}\label{luvv-2}
Let $v=|\u|$. Then
\begin{equation}
 (\cal{L}(u)v,v)=\int_\M |\u|^2\big(|\GV|^2+\Ric(V,V)\big)\,d\mu\ .
\end{equation}
\end{lem}
\begin{proof}
The function $u$ is of class $C^3$, hence $V$ is $C^2$.
In the open set $U$ we have $V(\Delta u+f(u))=0$, thus
\begin{equation}
\Delta(Vu)+f'(u)(Vu)=\Delta (Vu)-V(\Delta u)\ .
\end{equation}
Applying Lemma \ref{ident-V} directly to the righthand side of
(16) we get
\begin{equation}\label{LVu}
\Delta(Vu)+f'(u)(Vu)=\div(\GVS\u)+\inp{H_u}{\GV}+\Ric(\u,V)\ .
\end{equation}
The covariant derivative of $V$ is given by
\begin{equation}\label{GV-eq}
   \GV=\frac1{|\u|}H_u-\u\otimes\frac{(H_u\u)^{\ast}}{|\u|^3}\ .
\end{equation}
A computation shows that $\GV$ is orthogonal to the tensor
$\u\otimes\frac{(H_u\u)^{\ast}}{|\u|^3}$. Recalling that $v=|\u|=Vu$
we obtain
\begin{equation}\label{prodhugv}
  \inp{H_u}{\GV}=|\u|\inp{\frac1{|\u|}H_u-\u\otimes
     \frac{(H_u\u)^{\ast}}{|\u|^3}}{\GV}=v\,|\GV|^2\ .
\end{equation}
Let $W$ be any vector in the tangent space over a point of $U$.
Since $V$ is unitary we have
\begin{equation}
   \inp{\GVS V}{W}=\inp{V}{\nb_WV}=\frac12 W|V|^2=0\ .
\end{equation}
Thus $\div(\GVS\u)=\div(v\GVS V)$ vanishes identically. With the
help of \eqref{prodhugv} equation \eqref{LVu} turns into
\begin{equation}
  \Delta v+f'(u)v=v\,|\GV|^2+v\,\Ric(V,V)\ .
\end{equation}
Notice that $v$ vanishes in $\M-U$. Looking at the left-hand side
of the above identity as a distribution it becomes clear that
its support is contained in $\ov{U}$. Therefore, applying it on
$v\in H^1(\M)$ one obtains
\begin{equation}\label{Luvv}
 (\cal{L}(u)v,v)=\int_\M |\u|^2\big(|\GV|^2+\Ric(V,V)\big)\,d\mu\ ,
\end{equation}
which proves the Lemma.
\end{proof}

%
%

\begin{rem}
Let $p\in M$ be a non-critical point of $u$. The level set
$S=\{x\,|\,u(x)=u(p)\}$ is a regular hypersurface near $p$. It can
be seen that $\GV=A+(\nb_VV)\otimes V^\ast$, where $A:TS\to TS$ is
the \emph{shape operator} respect to $V$ of the second fundamental
form of the inclusion $S\subset M$. By setting  $c=|\nb_VV|$ the
squared norm of $\GV$ becomes
\begin{equation}
  |\GV|^2=|A|^2+c^2\ .
\end{equation}
Therefore $|\GV|^2$ is the sum of the square of the principal
curvatures of $S$ plus the square of the curvature of the flow of
$\u$.
\end{rem}
\begin{rem}
In the unidimensional case $\M=S^1$ a direct proof of instability
can be given. Endow $S^1$ with a metric so that $|S^1|=l$. Functions
on $S^1$ are identified with functions on $[0,l]$ satisfying certain
boundary conditions. In this case the Euler-Lagrange equation for
$\mathcal{E}$ is
\begin{equation}\label{auto-eq-lin}
  \left\{ \begin{aligned} &u''(t)+f(u(t))=0,\qquad\ 0<t<l \\
                           &u(0)=u(l) \   \\
                           &u'(0)=u'(l)\
\end{aligned}\right.
\end{equation}
Its linearization becomes $\cal{L}(u)v=v''+f'(u)v$.
Assume by contradiction that $u$ is a non-constant
local minimizer of $\Ee$. Then $(\cal{L}(u)v,v)\leq0$, and
due to Lemma \ref{luvv-2} we get $\cal{L}(u)v=0$.
Hence $v=|u'|$ is an eigenfunction associated to the zero eigenvalue.

A direct computation shows that $u'$ is also an eigenfunction of the
zero eigenvalue of $\cal{L}(u)$. Then $w=u'+|u'|$ is an
eigenfunction and since $w$ vanishes in an open interval the Unique
Continuation Theorem gives us $w\equiv 0$. Hence $u'\equiv 0$, what
goes against the hypothesis. This shows that the first eigenvalue of
$\cal{L}(u)$ is positive and there are no non-constant local
minimizers of $\Ee$.
\end{rem}

In view of Lemma \ref{luvv-2} the proof of Theorem \ref{tunstable}
is now immediate if we strengthen the hypothesis to $\Ric>0$ on
$\M$. Indeed, one can show that $\Ric>0$ on some open set of $\M$
suffices for the positivity of $(\cal{L}(u)v,v)$, by using the
Unique Continuation Theorem together with the contradiction
assumption that the first eigenvalue of $\cal{L}(u)$ is zero.

We will rather give a unified proof for the case $\Ric\geq 0$. This
requires a few more lemmas dealing with the more delicate case
$\GV=0$ and $\Ric=0$ on $U$. It will follow after a series of steps
rich on tricky details. The main ingredients are the level sets of
$u$ and the behaviour of the geodesics of $\M$ respect to the
critical points of $u$.

The remaining results of this section do not demand that $u$ be
bounded or belong to any particular Sobolev Space. We will skip
for a while any functional analytic concerns, and assume that
$\M$ is an arbitrary complete, not necessarily compact, Riemann
manifold, and $u$ is a \emph{classical} solution to equation
\eqref{autoeq}. The compactness of $\M$ will be implicitly invoked back
only in the proofs of Theorems \ref{tunstable} and \ref{tsemistable}.

For the next six Lemmas and Corollaries we thus assume
\begin{equation}
 |\GV|=0\qquad \text{in}\ U\ ,
\end{equation}
unless otherwise stated. In particular we obtain that
$V$ is a parallel vector field over $U$. From equation \eqref{GV-eq}
we also get
\begin{equation}\label{forma-Hu}
  H_u=V\otimes(H_u\,V)^\ast=\Delta u\,V\otimes V^\ast\ .
\end{equation}
For any $p\in\M$ define $\N_p$ as the component
of the level set $\{x\in\M\,|\,u(x)=u(p)\}$ that contains $p$.
\begin{lem}\label{u-level}
If $p\in U$ then $\N_p\subset U$. Further, $\N_p$ is a complete
geodesic riemannian submanifold of codimension 1 of $\M$ and
$|\u|>0$ is constant on $\N_p$.
\end{lem}
\begin{proof}
Let $U_p$ be a component of $U$ and $C_p$ a component of
$U_p\cap\N_p$ so that $p\in C_p$.
Clearly $C_p$ is a codimension 1 submanifold
of $\M$. If $X,Y\in T(C_p)\subset TU$ we have
\begin{equation}
   \inp{\nb_XY}{V}=X\inp YV-\inp Y{\nb_XV}=0
\end{equation}
for $V$ is parallel and $X,Y$ are orthogonal to $V$. This
shows that $C_p$ is geodesic.

Letting $q\in C_p$ and $X\in T_q(C_p)$, we have $\nb_X\u=H_u(X)=0$.
Therefore $\u$ is parallel and $|\u|\neq0$ is constant along $C_p$.
If $\ov{q}$ is an adherent point of $C_p$ then $\u(\ov{q})$ is
non-zero so that $\ov{q}\in U_p$. This shows that $C_p$ is closed in
$\M$, and since $U_p$ is open, $C_p$ is also open as a topological
subspace of $\N_p$. Therefore by the conexity we have
$C_p=\N_p\subset U_p$.

The geodesic completeness of $\N_p$ follows from the Theorem
of Rinow and Hopf \cite{dC} and the fact that $\M$ is complete.
\end{proof}

\begin{lem}\label{geod-nconst}
Let $\gm:\R\to\M$ be an arclength parametrized geodesic, and
$h(t)=u(\gm(t))$ for all $t\in\R$. Assume that $h$ is non-constant,
and let $(a,b)$ be a component of $\gm^{-1}(U)$.
Then\\
(a) $h$ is strictly monotone in $(a,b)$.\\
(b) Assume $a\in\R$, and let $p=\gm(a)$. Then $p$ is a critical
point of $u$ and $H_u(p)\neq0$.\\
(c) Under the same hypothesis as (b) let $r=b-a\in\R\cup\{+\infty\}$.
Then $(a-r,a)$ is also a component of $\gm^{-1}(U)$.
Further, $h(t)$ is simmetric respect to $t=a$, i.e.,
$h(a-s)=h(a+s)$ for all $s\in\R$. \\
(d) Under the same hypothesis as (c), assume also $b\in\R$.
Then $h$ is periodic of period $2r$.
\end{lem}
\begin{proof}
For all $t\in\R$ we have $h'(t)=\inp{\u}{\gm'(t)}$. This justifies
the existence of the interval $(a,b)$, since $h$ is non-constant.
For all $t\in(a,b)$ we can write $h'(t)=|\u|\inp{V}{\gm'(t)}$. Both
of $V$ and $\gm'$ are parallel along $\gm$, hence
$\inp{V}{\gm'(t)}=k$ is a constant in $(a,b)$. We must have
$k\neq0$, otherwise the geodesic $\gm$ would be entirely contained
in $\N_{\gm(t_0)}$, for any $t_0\in(a,b)$, and $h$ would be
constant. Hence $k$ and $|\u|$ are non-zero in $(a,b)$ and part (a)
is proved.

We compute the second derivative of $h$ for any $t\in(a,b)$,
\begin{equation}
h''(t)=\frac{d}{dt}\inp{\u}{\gm'(t)}=\inp{H_u(\gm'(t))}{\gm'(t)}=
   \Delta u(\gm(t))k^2\ ,
\end{equation}
in view of equation \eqref{forma-Hu}. Then $h(t)$ is a
solution to the $2^{\rm nd}$ order equation
\begin{equation}\label{equa-h}
  h''+ k^2\,f(h)= 0
\end{equation}
on $(a,b)$. If $a\in\R$, $h$ satisfies the initial
condition $h(a)=u(p)$, $h'(a)=0$. By uniqueness of the Initial Value
Problem the constant function $t\mapsto u(p)$ is not a solution of
that problem, and therefore $u(p)$ is not a root of $f$.
Hence, $h''(a)=-k^2\,f(u(p))\neq0$, and $H_u(p)$ does not vanish.
This concludes part (b).

Due to $h''(a)\neq0$ there is a small left open neighborhood
of $a$ where $h'(t)\neq0$, and hence $\gm(t)\in U$ for $t<0$ small.
Therefore there is a component of $\gm^{-1}(U)$ of
the form $(c,a)$, for some $c\in(-\infty,a)$.
Let $J=(0,\min\{r,a-c\})$.

We define $h_-(s)=h(a-s)$ and $h_+(s)=h(a+s)$ for all $s\in\R$.
Then $h_-(0)=h_+(0)=h(a)$, $h_-'(0)=h_+'(0)=0$. Further, for $s\in J$
there are suitable constants $k_-,k_+$ that play
the role of $k$ on \eqref{equa-h}:
\begin{align*}
  h_-''+k_-^2\,f(h_-)&=0\ ,\\
  h_+''+k_+^2\,f(h_+)&=0\ .
\end{align*}
Again uniqueness for this problem will give us $h_-\equiv h_+$ as
long as we show that $k_-^2=k_+^2$.

Let $V_-(s)=V(\gm(a-s))$ and $V_+(s)=V(\gm(a+s))$ for all
$s>0$ small. Both of $V_-$ and $V_+$ can be continously
extended by parallel transport along $\gm$ to vectors $\tilde{V}_-$
and $\tilde{V}_+$, respectively, on $T_p\M$. We claim that the
(unitary) vectors $\tilde{V}_-$ and $\tilde{V}_+$ are colinear.
The (symmetric polinomials on the) eigenvalues of the
continuous symmetric tensor $H_u$ are continuous. The
special form of $H_u$ on $U$, given by equation \eqref{forma-Hu},
implies that for all small $s>0$, $H_u(\gm(a\pm s))$ has a zero eigenvalue
of multiplicity at least $N-1$, which is inherited by $H_u(p)$.
The remaining eigenvalue of $H_u(p)$, $\Delta u(p)$,
has to be non-zero (after part (b)) and simple. This
is an open condition, and the eigenspace associated to this
eigenvalue varies continuously, close to $p$. It is generated
by $V$ on $U$, therefore, we have $\tilde{V}_-=\pm\tilde{V}_+$.
Since
\begin{equation}\label{k+k-}\begin{aligned}
  k_-&=\lim_{s\to0^+}\inp{V_-(s)}{\gm'(a-s)}=\inp{\tilde{V}_-}{\gm'(a)}\\
  k_+&=\lim_{s\to0^+}\inp{V_+(s)}{\gm'(a+s)}=\inp{\tilde{V}_+}{\gm'(a)}\ ,
\end{aligned}\end{equation}
we get $|k_-|=|k_+|$, hence $h_-(s)=h_+(s)$ for $s\in J$.
Critical points of $h_-$ and $h_+$ happen together in this range
and correspond to intersections of $\gm(t)$ with the border of $U$.
Therefore $0<s\mapsto\gm(a-s)$ cannot leave $U$ before $s=r$, and
since the argument is symmetric, we conclude that $a-c=r$ and
$\gm^{-1}(U)$ contains $(a-r,a)$ as a component, which proves part (c).

Part (d) is now immediate. Clearly the symmetry of $h(t)$
holds respect to any critical point of $h$. If $r=b-a$ is finite
then we get $h(a+r+s)=h(a+r-s)=h(a-r+s)$ for any $0<s<r$.
In particular, an inductive argument shows that $\{a+mr\,|\,m\in\mathbb{Z}\}$
are all critical points of $h(t)$. The period of $h$ is $2r$ since
it intercalates increasing with deacreasing intervals between
consecutive critical points.
\end{proof}
\begin{rem}\label{g-n-rem}
From part (c) of the Lemma we have $h'(a+s)=-h'(a-s)$, and picking
$s>0$ small we obtain
\begin{equation}
  h'(a+s)=k_+|\u|_{\gm(a+s)}=-k_-|\u|_{\gm(a-s)}=-h'(a-s)\ .
\end{equation}
Therefore, $k_-=-k_+$ and $\tilde{V}_-=-\tilde{V}_+$.
\end{rem}
As a consequence of Lemma \ref{geod-nconst} we get $H_u(p)\neq 0$
and $\Delta u(p)\neq0$ for any critical point $p$ of $u$, since there
is a point $q\in\M$ with $u(q)\neq u(p)$ and a geodesic $\gm(t)$
joining $p$ to $q$.
Further, the set of critical points of $u$ is $\partial U=M-U$.

We are now ready to give the
\begin{proof}[Proof of Theorem \ref{tunstable}]
By Lemma \ref{luvv-2} along with the condition $\Ric\geq0$ we deduce
that $(\cal{L}(u)v,v)\geq 0$. We will show that this inequality is
strict, so $u$ cannot be a local minimum of $\Ee$. The case where
$\GV\neq 0$ is straightforward from the Lemma, so we assume in the
sequel that $\GV\equiv0$ on $U$.

Suppose by contradiction that the first eigenvalue of $\cal{L}(u)$
is non-positive. Then $(\cal{L}(u)v,v)=0$ and $v$ must be an
eigenfunction of $\cal{L}(u)$ associated to the zero eigenvalue.
Since $f'(u)v$ is continuous, standard elliptic regularity
applied to
\begin{equation}\label{eigen-eq}
   \Delta v+f'(u)\,v=0 \qquad\,\mbox{on}\ \M
\end{equation}
gives us $v\in C^2(\M)$. Computing the gradient of $v$ in $U$ we obtain
\begin{equation}
  \nb v=\nb |\u|=H_u(V)=\Delta u\,V\ .
\end{equation}
Let $p$ be a critical point of $u$ and $\gm(t)$ be a geodesic
satisfying the hypotheses on Lemma \ref{geod-nconst}, so that
$\gm(0)=p$. Following  the notation in the proof of the Lemma we
have, by part (b), that $\Delta u(p)\neq 0$. On the other hand,
Remark \ref{g-n-rem} gives us
\begin{equation}
 \lim_{t\to0^+}V_{\gm(t)}=-\lim_{t\to0^-}V_{\gm(t)}\neq 0\ .
\end{equation}
This shows that $\nb v$ is not even continuous at $p$, what
contradicts the $C^2$ regularity of $v$.
The only remedy is granting that the first
eigenvalue of $\cal{L}(u)$ is positive, which finishes
the proof of the Theorem.
\end{proof}

Notice that $V$ defines a line subundle of $T\M|_U$ that can be
extended over $\partial U$ by taking the only simple eigenspace
of $H_u$ (associated to the non-zero eigenvalue) near critical
points. This justifies the next
\begin{cor}\label{line-bundle}
There exists a geodesic line bundle $\cal{I}\subset T\M$ so that
$\cal{I}|_{U}$ is spanned by $V$.
\end{cor}
Choose a point $p_0\in U$ and let $U_0$ be its correspondent
component of $U$. Denote $\N_0=\N_{p_0}$.
We would like to extend the field $V|_{U_0}$ to the whole of
$\M$ by means of the bundle $\cal{I}$.
The flow of such extension would, then, be generated by isometries,
and routine arguments would give us a covering map
$\varphi:\R\times\N_0\to\M$, from which one would quickly derive the
results of Theorem \ref{tsemistable}. This case has already been
researched in greater generality, for instance, in \cite{BN}.

Here is where the orientability of $\cal{I}$ comes in.
Clearly, such an extension of $V|_{U_0}$ is possible if and
only if $\cal{I}$ is orientable (as a real vector bundle).
Both of orientable and non-orientable cases can happen to $\cal{I}$,
leading to two different constructions for $\M$. In order
to keep generality and short the proofs, we give a definition
of $\varphi$ independent of $\cal{I}$.

For any $p\in\N_0$ let $t\in\R\mapsto \varphi_t(p)$ be the geodesic
defined by $\varphi_0(p)=p$ and $\varphi_0'(p)=V_p$. Then
$\varphi:\R\times\N_0\to\M$ is smooth.
\begin{lem}\label{isomet-U0}
There is an open interval $(a,b)$ so that
$\varphi:(a,b)\times\N_0\to U_0$ is an isometry.
\end{lem}
\begin{proof}
Let $(a,b)\ni0$ be the maximal interval for which $\varphi_t(p_0)$
belongs to $U_0$. If $q\in\N_0$ is any other point we see that
$u(\varphi_t(p_0))=u(\varphi_t(q))$ for $t\in\R$, since both functions
satisfy the same differential equation \eqref{equa-h} with same
initial conditions. Due to Lemma \ref{geod-nconst}
it follows that $(a,b)$ keeps the maximality
property above stated, for any $q\in\N_0$.

Since $V$ is parallel and equals $\varphi_0'(p)$ on $p$, it holds
$\varphi_t'(p)=V_{\varphi_t(p)}$ for all $t\in(a,b)$. Therefore
$t\mapsto\varphi_t$ are integral curves of $V|_{U_0}$. Two such
curves do not intersect, and because $u(\varphi_t(p))$ is monotone
the curve $\varphi_s(q)$ cannot be a reparametrization of
$\varphi_t(p)$, for any $(s,q)\in(a,b)\times\N_0$ with $q\neq p$.
This concludes injectivity of $\varphi:(a,b)\times\N_0\to U_0$.
Notice that $\varphi$ is the flow of $V$ restricted to $\N_0$,
hence it is an isometry with its image.
The set $\varphi((a,b)\times\N_0)$ is open.

Now we show that the image of $\varphi$ is closed in $U_0$.
Let $q\in U_0$ be an adherent point of $\varphi((a,b)\times\N_0)$,
and $\sg:[0,1]\to U_0$ be a smooth curve
with $\sg(0)=p_0$, $\sg(1)=q$. Let $I=\sg^{-1}(U_0)$, $I$ is
open in $[0,1]$ and non-empty.
Using that $\varphi$ is a local isometric coordinate chart
one see that $I$ is closed, hence $I=[0,1]$ and $q$
belongs to the image of $\varphi$.
The image of $\varphi$ is then open and closed in $U_0$, and by
conexity, we have $\varphi((a,b)\times\N_0)=U_0$.
\end{proof}

Following the notation of Lemmas \ref{geod-nconst} and
\ref{isomet-U0} we consider the case $b\in\R$.
Then $\varphi_b(\N_0)\subset\partial U_0$. Let $\hat p=\varphi_b(p_0)$.
We have $\varphi_b(\N_0)=\N_{\hat p}$, since $\varphi_t$ preserves
level sets of $u$. Surprisingly, $\N_{\hat p}$ may not be isometric
to $\N_0$. This question relates to whether the curve
$t\mapsto\varphi_t(p)$ does \emph{leave} $U_0$ when it crosses
the border at $t=b$.

Let $U_1$ be the component of $U$ that contains
$\varphi_t(\N_0)$ for all $b<t<2b-a$.
\begin{lem}\label{U0-border}
$\N_{\hat p}$ is a geodesic complete submanifold of $\M$.
The map $\varphi_b:\N_0\to\N_{\hat p}$ is a local isometry.
It is a bijection if and only if $\cal{I}|_{U_0\cup\N_{\hat p}}$ is
orientable, and it holds $U_1\neq U_0$.
Otherwise $\varphi_b$ is a two-fold covering map
onto $\N_{\hat p}$ and $U_1=U_0$.
\end{lem}
\begin{proof}
Let $p\in\N_0$ and $\cal{V}\ni\varphi_b(p)$ be a simply connected
open neighborhood of $\varphi_b(p)$. There is a local trivialization
of $\cal{I}|_{\cal{V}}$ by means of a unitary parallel vector field
$\tV$, so that $\tV_{\varphi_b(p)}=\varphi_b'(p)$. By continuity,
$\varphi_b'(q)=\tV_{\varphi_b(q)}$ for any $q\in\varphi_b^{-1}(\cal{V})$.
Again, uniqueness of the parallel trasport along a curve subject
to the same initial conditions gives us
$\tV_{\varphi_{b+s}}(q)=\varphi_{b+s}'(q)$ for all $s$ small enough.
Restricting $\cal{V}$ if necessary we see that $\varphi$ is the flow
of a unitary killing field defined on the open set
$\varphi((a,b+\varepsilon)\times\varphi_b^{-1}(\cal{V}))\cup\cal{V}$,
for some $\varepsilon>0$ small. Hence $\varphi_b$ is a local
isometry of $\N_0$ onto $\N_{\hat p}$. From that it also follows
that $\N_{\hat p}$ is geodesic and complete.

Now assume $\varphi_b$ is injective. Then $(t,q)\in(a,b]\times\N_0
\mapsto\varphi_t'(q)$ is a well defined trivialization of
$\cal{I}|_{U_0\cup\N_{\hat p}}$, so it is orientable.
If $\varphi_t(p)$ belongs to $U_0$ for some $t\in(b,2b-a)$ then
there is $s\in(a,b)$ and $q\in\N_0$ with $\varphi_s(q)=\varphi_t(p)$.
Both geodesics have velocities on the bundle $\cal{I}$,
so they must be opposite since $u(\varphi_t(p))$ is decreasing on $t$.
Therefore $\varphi_t(p)$
is a backward reparametrization of $\varphi_s(q)$ and we get
$\varphi_b(p)=\varphi_b(q)$, contradicting injectivity.
Hence there must be $U_0\neq U_1$.

On the other hand, if there are distinct points
$p,q\in\N_0$ with $\varphi_b(p)=\varphi_b(q)$ one clearly has
$\varphi_b'(p)=-\varphi_b'(q)$, since both velocities lie in the
same fiber of $\cal{I}$ and cannot be equal. Therefore no
orientation of $\cal{I}|_{U_0}$ can be extended to a larger
set on $\M$ containing $\N_{\hat p}$, i.e., $\cal{I}|_{U_0\cup\N_{\hat p}}$
is non-orientable. In this case it holds $\varphi_{2b}(p)=q$, hence
$\varphi_{2b}(\N_0)=\N_0$, what indicates that $U_0=U_1$.
Restricting $\varphi_b$ to suitable vicinities $\cal{V}_p$, $\cal{V}_q$
of $p$ and $q$, respectively, we may write
$\varphi_{2b}|_{\cal{V}_p}=(\varphi_b|_{\cal{V}_q})^{-1}
\circ\varphi_b|_{\cal{V}_p}$,
what shows that $\varphi_{2b}$ is locally an isometry
\emph{without fixed points} and $\varphi_{2b}^2=Id_{\N_0}$.
This finishes the proof that $\varphi_b:\N_0\to\N_{\hat p}$ is
a two-fold covering map.
\end{proof}
Recall that an involution of a riemannian manifold is an isometry
$I$ such that $I^2=id$.
\begin{lem}\label{regular-cover}
$\varphi:\R\times\N_0\to\M$ is a regular isometric covering map.
Denote by $K=Aut(\R\times\N_0,\varphi)$ the group of covering
transformations of $\varphi$. Then, if $\cal{I}$  is orientable,
$K$ is either trivial or ciclic generated by the metric product
of a translation of $\R$ with an isometry of $\N_0$. If $\cal{I}$ is
not orientable $K$ is generated by at most two involutions of
$\R\times\N_0$.
\end{lem}
\begin{proof}
If $u$ has no critical points then $U_0=U=\M$ and $\varphi$ is
the (regular) trivial covering map, $\cal{I}$ is orientable and $K=\{Id\}$.
Otherwise $\partial U_0\neq\emptyset$ and
we assume $b$ on Lemma \ref{isomet-U0} is finite.

Following  Lemma \ref{U0-border} we let $\N_{\hat p}=\varphi_b(\N_0)$ be
a component of the border of $U_0$. If there is another
component $U_1$ of $U$ that cobounds $U_0$ through $\N_{\hat p}$
then we can choose $p_1\in U_1$ with $u(p_1)=u(p_0)$ and
let $\N_1=\N_{p_1}$. Let $\psi:(a,b)\times\N_1\to U_1$ be the map
analogous to $\varphi$. It can be seen from the proof of Lemma
\ref{U0-border} that $\varphi_t'(p)\in\cal{I}_{\varphi_t(p)}$
for all $t\in\R$, $p\in\N_0$.
Then $\psi_b(\N_1)=\varphi_b(\N_0)=\N_{\hat p}$. It is clear
that $\varphi_{2b}(\N_0)=\N_1$ and
$\varphi_{b+s}(p)=\psi_{b-s}(\varphi_{2b}(p))$ for all $s\in\R$,
$p\in\N_0$.
Therefore $\varphi$ is an isometry from $(a,2b-a)\times\N_0$
onto $U_0\cup\N_{\hat p}\cup U_1$.

On the other hand, if $U_0$ \emph{self-bounds} at $\N_{\hat p}$
as described by Lemma \ref{U0-border}, the function $\psi$ above
defined equals $\varphi$, and $\N_1=\N_0$.
Hence $\varphi:(a,2b-a)\times\N_0\to U_0\cup\N_{\hat p}$ is a
two-fold isometric covering map.

If $a=-\infty$ we are done. Otherwise there is another component
$\N_{\hat q}$ of $\partial U_1$, $\N_{\hat q}\neq\N_{\hat p}$.
The above constructions can be repeated, extending the isometric
covering property of $\varphi$ to the interval $(a,3b-2a)$.
This can also be performed backwards on $t$, starting on $t=a$.
An inductive argument gives us that $\varphi:\R\times\N_0\to\M$
is a covering map, and a local isometry.

If $\varphi$ is injective we have again the trivial covering, and
$K=\{Id\}$. In this case one clearly has $\cal{I}$ orientable.
We assume in the remaining of this proof that $\varphi$ is not injective.

Suppose first that $\cal{I}$ is orientable.
Let $\varphi_{t_1}(p_1)=\varphi_{t_2}(p_2)$ for some
$(t_1,p_1),(t_2,p_2)\in\R\times\N_0$ distinct.
Then $\varphi_{t_1}'(p_1)=\varphi_{t_2}'(p_2)$, so $\varphi_t(p_1)$
is an orientation preserving reparametrization of $\varphi_s(p_2)$.
There is $\tau>0$ with $\varphi_\tau(\N_0)=\N_0$, and $\tau$ can
be taken the smallest positive number with such property.
Then $\varphi_\tau$ is an isometry of $\N_0$.

Consider the automorphism of the covering space $\R\times\N_0$
given by $g_\tau(t,p)=(t-\tau,\varphi_{\tau}(p))$. A quick computation
shows that the subgroup generated by $g_\tau$ acts transitively
on the preimage $\varphi^{-1}(q)$ for all $q\in\M$.
Since $K$ is completely defined by some subgroup of the permutations
group of $\varphi^{-1}(q)$ it becomes $K=\{g_\tau^n\,|\,n\in\mathbb{Z}\}$,
and the covering map is regular.

Now consider $\cal{I}$ not orientable. Reasoning similarly to
the previous case we can find $C\neq0$ so that
$\varphi_C:\N_0\to\N_{\hat p}$, $\hat p=\varphi_C(p_0)$,
is a two-fold covering, and $\varphi_{2C}:\N_0\to\N_0$ is an involution.
We can pick $C$ so that $|C|>0$ is minimum.
Then $g_C(t,p)=(2C-t,\varphi_{2C}(p))$ is an involution of $\R\times\N_0$
and a covering transformation.
If $\varphi$ is a two-fold covering then the orbits of $\{Id,g_C\}$
acting on $\R\times\N_0$ are all the preimages of points of $\M$.
Hence $\varphi$ is regular and $K=\{Id,g_C\}$.

If $\varphi$ is not a two-fold covering let $(t_2,p_2)$, $(t_1,p_1)$ and
$g_C(t_1,p_1)$ be three distinct points in the preimage of a fixed
point $q\in\M$. The velocities of the geodesics $s\mapsto\varphi_s(p_1)$
and $s\mapsto\varphi_{2C-s}(\varphi_{2C}(p_1))$ are opposite over
$q$, and we can assume, without loss of generality, that
$\varphi'_{t_2}(p_2)=\varphi'_{t_1}(p_1)$.
Again there is $\tau>0$ minimum such that $\varphi_\tau(\N_0)=\N_0$
and $\varphi'_\tau(p)=V_{\varphi_\tau(p)}$ for any $p\in\N_0$.
Define $g_\tau$ as in the $\cal{I}$ orientable case.

Now let $(t,p)$ be any point in $\varphi^{-1}(q)\ni (t_1,p_1)$.
If $\varphi'_t(p)=\varphi'_{t_1}(p_1)$ then there is an integer
$n$ such that $(t,p)=g_\tau^n(t_1,p_1)$.
Otherwise $(t,p)=g_\tau^n\circ g_C(t_1,p_1)$. This shows that the
action of $K$ is transitive on the preimages and the covering map
is regular. Further $K$ is generated by $\{g_\tau,g_C\}$. A careful
check traveling forth and back on the geodesics $t\mapsto\varphi_t(p)$
reveals that
$\varphi_\tau\circ\varphi_{2C}\circ\varphi_\tau\circ\varphi_{2C}=Id_{\N_0}$.
Defining $D=C-\frac\tau2$ and $g_D(t,p)=(2D-t,\varphi_{2D}(p))$
we see that $g_D=g_\tau\circ g_C$ is an involution of $\R\times\N_0$
and $\{g_C,g_D\}$ generates $K$. This finishes the proof of the Lemma.
\end{proof}

\begin{proof}[Proof of Theorem \ref{tsemistable}]
Let $u$ be a non-constant critical point of $\Ee$ with\\
$(\cal{L}(u)v,v)=-(\Ee''(u)v,v)=0$.
Clearly the manifold $\N$ in the Theorem stands for $\N_0$.

The proof then follows from the sequence of the Lemmas and
Corollaries numbering from \ref{u-level} through
\ref{regular-cover}. The assertion $\M\simeq(\R\times\N)/K$ is a
standard fact in Topology $\cite{Ms}$ and the metric is induced from
$\R\times\N$ through the local isometry $\varphi$.
\end{proof}

\section{Existence of nonconstant minimizers}

This section is devoted to show that if $\M$ fails to have
non-negative Ricci curvature then Theorem \ref{tunstable} may not
hold. This will be accomplished by showing  that there are
non-convex surfaces for which $\Ee_{\varepsilon}$  has non-constant
local minimizers, for $\varepsilon$ small enough.

The procedure we follow consists of finding the limit of the
energies $\Ee_\varepsilon$ in the sense of $\Gamma-$convergence and
then using a result of De Giorgi which roughly states that close (in
some specified topology) to an isolated minimizer of the
$\Gamma$-limit problem there is a minimizer of the original one.

Throughout this section, $\M$ will denote a surface diffeomorphic to
$S^2$. For the reader's convenience we give the definition of
$\Gamma-$convergence which is going to be used.

A family $\{\Lambda_{\varepsilon}\}_{0<\varepsilon\leq\varepsilon_0}$ of
real-extended functionals defined in $L^1(\M)$ is said to
$\Gamma$-converge in $L^1(\M)$ , as $\varepsilon\to 0$,
to a functional $\Lambda_0:L^1(\M)\longrightarrow \R \cup
\{\infty\}$, if:
\begin{itemize}
\item For each $v\in L^1(\M)$ and for any family
 $\{ v_\varepsilon\}$ in $L^1(\M)$ such that $v_\varepsilon\to v$ in
$L^1(\M)$, as $\varepsilon\to 0$, it holds that
$\Lambda_0(v)\leq\lim\inf_{\vspace*{-6mm}\varepsilon\to 0}
\Lambda_{\varepsilon}(v_\varepsilon)$.
\item For each $v\in L^1(\M)$ there is a family
$\{ w_\varepsilon\}$ in $L^1(\M)$ such that $w_\varepsilon\to v$ in
$L^1(\M)$, as $\varepsilon\to 0$ and  $\;\;\Lambda_0(v)\geq
\lim\sup_{\vspace*{-6mm}\varepsilon\to 0}\Lambda_{\varepsilon}
(w_\varepsilon)$.
\end{itemize}

Convergence in this sense will be denoted by
${\Gamma}^-\lim_{\varepsilon \rightarrow
0^+}\Lambda_{\varepsilon}=\Lambda_0$. The definitions and  results
we need about functions of bounded variation defined on $\M$ are
provided below.

We set
\begin{equation}
\cal{G}(\M)\stackrel{\rm def}{=}\{g\ |\ g\ \text{is a}\ C^1\
\text{section of}\ T\M,\  |g(x)|\leq 1,\ \forall\ x\in\M\}
\end{equation}
and let  $\cal{H}^{N}$ denote the usual $N$-dimensional Hausdorff
measure.

Given  $u:\M  \rightarrow \R$  we define
\begin{equation}
 |D u|(\M)\stackrel{\rm def}{=}\sup_{g \in \cal G(\M)
}\int_\M u\,\div(g)\,d\cal{H}^2\ .
\end{equation}
A real function  $u \in L^1(\M)$ has bounded variation
in $\M$ if $
 |D u|(\M ) < \infty
.$ See \cite{Giu} when $\M$ is a bounded domain in $\R^N$. The set
$$ BV(\M ) \stackrel{\rm def}{=} \{ u : \M \rightarrow \mathbb{R} ;
\text{ } u \in L^1 (\M) \text{ and}\, |D u| (\M ) < \infty \}$$ is a
Banach space with the norm $\| u \|_{BV} =
   \| u \|_{L^1} + |D u | (\M ) $.

Letting ${\chi}_A$ denoting the characteristic function of a set
$A \subset \M$ we have
\begin{equation}
|D {\chi}_A|(\M )=\sup_{g\in\cal G(\M)}\int_A\div(g)\,d\cal{H}^2\ .
\end{equation}
The perimeter of a set $A \subset \M $ is defined by $\text{Per}_\M
(A):= |D {\chi}_A|(\M )$. If the border of $A$ in $\M$ is at least
$C^2$ then $|D\chi_A |(\M)=\mathcal{H}^1(\pa A\cap\M)$.

Throughout this section  we assume that the potential $F$ in (3)
satisfies:
\begin{itemize}
  \item $F: {\R} \rightarrow \mathbb{R}$ is $C^2$
  \item $F \geq 0$ and $F(t) = 0 $ if and only if $t
         \in \{ \alpha,\beta \}$,\, $\alpha<\beta.$
  \item $\exists\,\, t_0>0,\;c_1>0,\;c_2>0,\;k>2$ such that $c_1t^k \leq F(t) \leq
c_2t^k,$ for $|t| \geq t_0.$
\end{itemize}

For convenience we denote the space of functions of bounded
variation in $\M$ taking only two values, $\alpha$ and $\beta$, by
$\bvab$.

The computation of the $\Gamma-$limit of $\Ee_{\varepsilon}$
when $\M$ is a bounded domain in $R^N$ is standard by now.
However no such result is available in the literature when $\M$
is a surface.
Nevertheless the proof  found in \cite{AD} can be adapted
to our case in a natural manner thus yielding
\begin{thm}\label{gam-conv}
 Let  $\Ee_{\varepsilon} :
  L^1 (\M) \rightarrow\R$ be defined by
\begin{equation}
 \Ee_{\varepsilon}(u) = \left\{\begin{array}{ll}
\displaystyle\int_\M \left[\varepsilon\frac{|\nabla u |^2}2
    -\varepsilon^{-1}F(u)\right]d\mathcal{H}^{2}
    &\text{if } u \in H^1(\M) \\ \infty
    &\text{if } u \in L^1(\M)\backslash H^1(\M)
\end{array} \right.
\end{equation}
Then $ {\Gamma}^- \lim_{\varepsilon \rightarrow 0^+}
\Ee_{\varepsilon}=\Ee_{0}$ where
\begin{equation}
\Ee_{0}(u)= \left\{
\begin{array}{ll} \lambda\,|D{\chi}_{\{u=\alpha\}}|(\M) & \text{if
} u  \in\bvab
\\ \infty & \text{otherwise}
\end{array} \right.
\end{equation}
and
\begin{equation}\label{constant-s}
\lambda= \displaystyle\int_0^1 \sqrt{F(s)}\,ds\ .
\end{equation}
\end{thm}

We say that $v_0\in L^1(\M)$ is an $L^1$-local minimizer
of the functional $\Lambda_0:L^1(\M)\longmapsto\R\cup\{\infty\}$
if there is $r>0$ such that
\[ \Lambda_0(v_0)\leq \Lambda_0(v) \ \ \mbox{whenever} \ \
   0 < \| v-v_0\|_{L^1(\M)} <r\ .\]
Moreover if $\Lambda_0(v_0)<\Lambda_0(v)$ for
$0<\|v-v_0\|_{L^1(\M)}<r$, then $v_0$ is called an
isolated $L^1$-local minimiser of $\Lambda_0$.

The next result, which we use in order to find a family of
minimizers for \eqref{variation-eps}, is due to De Giorgi and can be
found in its abstract form in \cite{G}. A proof, with the hypotheses
on $F$ given above, can be found in \cite{KS}, since the replacement
of Lebesgue measure with Haussdorf measure does not affect the
arguments used.

\begin{thm}\label{deG-t}
Suppose that a sequence of real-extended functionals
$\{\Lambda_\varepsilon\}$ and $\Lambda_0$ satisfy

\begin{description}
\item[(\rm i)]
 ${\Gamma}^-
\lim_{\varepsilon \rightarrow 0^+} \Lambda_{\varepsilon}
=\Lambda_0$
 \item[(\rm ii)]Any sequence $\{ v_\varepsilon\}_{\varepsilon >0}$ such that
$\Lambda_\varepsilon(v_\varepsilon)\leq C <\infty$ for all
$\varepsilon
>0$, is compact in $L^1(\M)$.
\item[\rm (iii)] There exists an isolated $L^1$-local minimizer
$v_0$ of $\Lambda_0$\,.
\end{description}

Then  $\exists\; \varepsilon_0>0$ and a family $\{
v_\varepsilon\}_{0 < \varepsilon\leq\varepsilon_0}$ such that
\begin{itemize}
\item  $v_\varepsilon$ is an $L^1$-local minimiser of
$\Lambda_\varepsilon$, and
\item  $\|v_\varepsilon-v_0\|_{L^1(\M)}\to 0$\,, \ as \
$\varepsilon\to 0$.
\end{itemize}
\end{thm}

The growth condition on $F$ is required in order to have the
hypothesis on compactness (ii) satisfied. We also take, without loss
of generality,  $\lambda=1$ on equation \eqref{constant-s}.

For any $u\in\bvab$ we denote by $\gm$ its boundary curve, i.e.,
$\gm=\partial\{p\in\M\,|\,u(p)=\al\}$. Similarly, for
any such $\gm$ there are exactly two distinct functions in $\bvab$
with $\gm$ as boundary curve. It holds $\Ez(u)=|\gm|$. Given $r>0$
there exists $\tilde{u}\in\bvab$ so that $\tilde{\gm}$ is the
disjoint union of a finite number of smooth closed curves satisfying
\begin{itemize}
\item $\|u-\tilde{u}\|_{BV}<r$;
\item $|\gm|\geq |\tilde{\gm}|$.
\end{itemize}
We set
\begin{multline}
\bvabs=\{u\in\bvab\,|\,\gm\subset\M\ \\
\text{is a smooth 1-dimensional submanifold}\}\ .
\end{multline}

Now we assume that a simple closed geodesic $\gz$ is separable,
i.e., $\M-\{\gz\}$ has two components. Let $u_0\in\bvabs$
 be the function associated to $\gz$ so that $u_0=\al\chi_{M_{\al}}+\be\chi_{\M_{\be}}$ with
$M_i=\{p\in\M\,|\,u_0(p)=i\}$\, $(i=\al,\be)$.
\begin{thm}\label{min-geod}
Under the hypotheses and notation of Theorem \ref{tstable} it
holds that $u_0$ is an $L^1(\M)$-local isolated minimizer of
$\Ee_0.$
\end{thm}
\begin{proof}
Let $\vv$ be the neighborhood constructed in preparation for
Lemma \ref{length-K}. We choose $0<\delta_0<\delta$ and
define $\vz=\varphi([-\delta_0,\delta_0]\times\gz)$.
We claim that any $r>0$ with
\begin{equation}\label{lisupeps}
  r<|\be-\al|\,\delta_0\,\min\left\{\delta-\delta_0,
          \frac{|\gz|}2\right\}
\end{equation}
will verify $\Ez(u)>\Ez(u_0)$ whenever $u\in\bvab$ and
$0<\|u-u_0\|_{L^1}<r$.

The discussion prior to the theorem allows us to restrict our
attention to competing functions $u\in\bvabs$. Let $\gm$ be the
boundary curve of a given $u$. A differential topology argument (see
\cite{GP}) allows us to consider $\gm$ in generic position with $\pa
\vz$ and $\pa \vv$, or equivalently, $\gm$ is transversal to the
boundaries of $\vz$ and $\vv$. In particular, each connected
component of $\gm\cap \vz$ is diffeomorphic to either $S^1\subset
\mbox{int}\,\vz$ or $[0,1]\subset \vz$ and endpoints contained in
$\pa \vz$. We define
\begin{equation}\begin{aligned}
 D&=\{\sg\ |\ \sg\ \text{is a component of}\ \gm\cap \vz\}\ ,\\
 I&=\bigcup_{\sg\in D}\bsg\ \subset\ \gz\ .
\end{aligned}\end{equation}

\begin{lem}\label{auxiliar}
Let $u\in\bvabs$ with $\|u-u_0\|<r$. Then
$|I|>\max\left\{|\gz|-(\delta-\delta_0),\dfrac{|\gz|}2\right\}$.
\end{lem}
\begin{proof}
For each $\sg\in D$, $\bsg$ is a closed segment of $\gz$.
Hence,
\begin{equation}
 J\stackrel{\rm def}{=}\gz-I=\bigcup_{i=1}^m J_i\ ,
\end{equation}
where each $J_i$ is an open interval of $\gz$, and the
$J_i$'s are pairwise disjoint. The construction leading to $J$
clearly yields
\begin{equation}
 \gm\cap\varphi([-\delta_0,\delta_0]\times J_i)=\emptyset\qquad
\text{for}\ 1\leq i\leq m\ .
\end{equation}
Therefore, $u$ is constant in $\varphi([-\delta_0,\delta_0]\times J_i)$.
Since $u_0$ switchs its value over $J_i$ we conclude that
$|u-u_0|=|\be-\al|$ in one of the regions $\varphi([-\delta_0,0]\times J_i)$
or $\varphi([0,\delta_0]\times J_i)$. Applying Lemma \ref{length-K}
part (b) we derive
\begin{equation}
  \|u-u_0\|_{L^1(\varphi([-\delta_0,\delta_0]\times J_i))}>
    |\be-\al|\delta_0\,|J_i|\ .
\end{equation}
Thus
\begin{align}
   r>\|u-u_0\|_{L^1}&>\sum_{i=1}^m|\be-\al|\delta_0\,|J_i|=
     |\be-\al|\delta_0(|\gz|-|I|)\\
   &\Rightarrow\ |I|>|\gz|-\frac{r}{|\be-\al|\delta_0}\ .
\end{align}
Together with \eqref{lisupeps} the above inequality  readily
implies the Lemma.
\end{proof}
We set a little more notation: for any $\sg\in D$ let $\rho=\rho(\sg)$
be the component of $\gm$ that contains $\sg$ as an arc.
We are led to three cases:\\
(i) If there is some $\rho(\sg)\not\subset \vv$ then there is an arc
$\tilde{\sg}\subset\rho$ joining a point of $\pa \vz$ to a point of
$\pa \vv$. Lemma \ref{length-K} (part (a)) gives us
$|\tilde{\sg}|\geq \delta-\delta_0$ and then
\begin{equation}\begin{aligned}
  |\gm|\geq |\tilde{\sg}|&+\sum_{\sg\in D}|\sg|\geq\delta-\delta_0+|I|\\
                  &> |\gz|\ ,
\end{aligned}\end{equation}
in view of Lemma \ref{auxiliar}.\\
(ii) If there is some $\rho(\sg)\subset \vv$ that is freely
homotopic to $\gz$ within $\vv$ then the intersection number
of $\rho$ with any geodesic ray $t\mapsto\varphi_t(x)$ is $\pm1$.
Denoting by $\brho$ the projection of $\rho$ over $\gz$
we get $\brho=\gz$.
Hence, Lemma \ref{length-K} part (a2) gives us
$|\gm|\geq|\rho|\geq|\gz|$. The strictness $|\gm|>|\gz|$
comes from $\|u-u_0\|_{L^1}>0$, since there must be
another component $\rho'\neq\rho$ of $\gm$ or $\rho$ is not
equal to $\gz$.\\
(iii) Assume that neither (i) nor (ii) occurs. If for some $\sg\in
D$ we have $\brho=\gz$ we conclude similarly to case (ii) above,
hence $|\gm|>|\gz|$. Otherwise, let $p$ and $q$ be points of $\rho$
so that their projections over $\gz$ are the end points of the
segment $\brho\subset\gz$. Let $\sg_1$ and $\sg_2$ be the two
distinct arcs of $\rho$ joining $p$ and $q$ ($\sg_i\subset \vv$,
$i=1,2$), with projections respectively $\bar{\sg_1}$ and
$\bar{\sg_2}$. Since the intersection number of $\rho$ with the ray
$t\mapsto\varphi_t(x)$ is 0 we have $\bar{\sg_1}=\bar{\sg_2}=\brho$.
Hence $|\rho|=|\sg_1|+|\sg_2|>2|\bar{\sg_1}|$. Fixing $\rho$ we see
that any $\sg\in D$ that is an arc of $\rho$ satisfies
$\bsg\subset\bar{\sg_1}$. Then,
\begin{equation}
  \left|\bigcup_{\sg\in D\,,\sg\subset\rho}
  \!\!\!\bsg\right|\leq |\bar{\sg_1}| <\frac12|\rho|\ ,
\end{equation}
from which we derive
\begin{equation}
   |\gm|\,\,=\!\!\!\!\!\!\!\sum_{\begin{array}{c}\rho\ \text{{\footnotesize
a component}}\\
       \text{{\footnotesize of}}\ \gm\end{array}}
\!\!\!\!\!\!\!\!\! |\rho|>2|I|>|\gz|\ .
\end{equation}
Therefore $\Ez(u)=|\gm|>|\gz|=\Ez(u_0)$ if $0<\|u-u_0\|_{L^1}<r$
and the theorem is proved.
\end{proof}

\begin{proof}[Proof of Theorem \ref{tstable}]
As mentioned before, Theorem \ref{tstable} is just an application of
Theorem \ref{deG-t} for $\Lambda_\varepsilon=\Ee_\varepsilon$, whose
hypotheses we now verify. Indeed (i) is nothing but Theorem
\ref{gam-conv} and (ii) may be found in \cite{S}, for instance.
Although the proof of (ii) in \cite{S} is rendered for $\M$ a
bounded domain in $\R^N$ the proof holds equally well in our case.

As for (iii) it has been verified in Theorem \ref{min-geod} above.
\end{proof}

The following result seems to be known, though we have not been able
to find it in the literature. It is a consequence of the procedure
used in this section along with Theorem \ref{tunstable}.
\begin{lem}
Let $\M$ be a compact Riemann surface with no boundary
and having nonnegative Gaussian curvature.
Then $\M$ has no closed nonintersecting isolated minimizing
geodesic.
\end{lem}

\end{document}